\documentclass[10pt]{article}
\begin{document}

\newcommand{\nc}{\newcommand}

\newtheorem{lemma}{Lemma}[section]
\newtheorem{theorem}[lemma]{Theorem}
\newtheorem{proposition}[lemma]{Proposition}
\newtheorem{corollary}[lemma]{Corollary}
\newtheorem{remark}[lemma]{Remark}
\newtheorem{example}[lemma]{Example}
\newtheorem{hypothesis}[lemma]{Hypothesis}
\newtheorem{notation}[lemma]{Notation}
\newtheorem{definition}[lemma]{Definition}
\newtheorem{conclusion}[lemma]{Conclusion}

\nc{\QED}{\mbox{}\hfill \raisebox{-0.2pt}{\rule{5.6pt}{6pt}\rule{0pt}{0pt}}
         \medskip\par}
\newenvironment{Proof}{\noindent
   \parindent=0pt\abovedisplayskip = 0.5\abovedisplayskip
   \belowdisplayskip=\abovedisplayskip{\bf Proof. }}{\QED}
\newenvironment{Proofof}[1]{\noindent
   \parindent=0pt\abovedisplayskip = 0.5\abovedisplayskip
   \belowdisplayskip=\abovedisplayskip{\bf Proof of #1. }}{\QED}
\newenvironment{Example}{\begin{example}
               \parindent=0pt \rm}{\QED\end{example}}
\newenvironment{Remark}{\begin{remark}
               \parindent=0pt\rm }{\QED\end{remark}}
\newenvironment{Definition}{\begin{definition}
               \parindent=0pt \rm}{\QED\end{definition}}
\newenvironment{Conclusion}{\begin{conclusion}
               \parindent=0pt \rm}{\QED\end{conclusion}}

\def\ba{\begin{array}}
\def\ea{\end{array}}

\def\be{\begin{equation}}
\def\ee{\end{equation}}
\def\vs5{\vspace{0.5cm}}
\def\lab{\label}
\def\bthm{\begin{theorem}}
\def\ethm{\end{theorem}}
\def\bhyp{\begin{hypothesis}}
\def\ehyp{\end{hypothesis}}
\def\bP{\begin{Proof}}
\def\eP{\end{Proof}}
\def\bPof{\begin{Proofof}}
\def\ePof{\end{Proofof}}
\def\brem{\begin{remark}}
\def\erem{\end{remark}}
\def\bex{\begin{example}}
\def\eex{\end{example}}
\def\bcor{\begin{corollary}}
\def\ecor{\end{corollary}}
\def\bd{\begin{definition}}
\def\ed{\end{definition}}
\def\bprop{\begin{proposition}}
\def\eprop{\end{proposition}}
\def\blem{\begin{lemma}}
\def\elem{\end{lemma}}
\def\beq{\begin{eqnarray}}
\def\eeq{\end{eqnarray}}
\def\beqs{\begin{eqnarray*}}
\def\eeqs{\end{eqnarray*}}

\def\a{\alpha}
\def\b{\beta}
\def\s{\sigma}
\def\Sig{\Sigma}
\def\d{\delta}
\def\l{\lambda}
\def\hs{\hat{\sigma}}
\def\hT{\hat{T}}
\def\hU{\hat{U}}
\def\eps{\epsilon}
\def\veps{\varepsilon}
\def\benum{\begin{enumerate}}
\def\eenum{\end{enumerate}}
\def\bit{\begin{itemize}}
\def\eit{\end{itemize}}
\def\la{\langle}
\def\ra{\rangle}
\def\wto{\rightharpoonup}
\def\bwto{\buildrel{w(\mu_h,\mu)}\over\longrightarrow}
\def\bYto{\buildrel{Y(\mu_h,\mu)}\over\longrightarrow}
\def\bsto{\buildrel{s(\mu_h,\mu)}\over\longrightarrow}
\def\bswto{\buildrel{sw(\mu_h,\mu)}\over\longrightarrow}

\def\Cinf{C^\infty}
\def\Linf{L^\infty}
\def\supp{{\rm supp}}
\def\setm{\setminus}
\def\id{{\rm id}}

\def\dis{\displaystyle}

\def\R{{\bf R}}
\def\N{{\rm N}}
\def\Ra{{\rm R}}
\def\M{{\cal M}}
\def\B{{\cal B}}
\def\A{{\cal A}}
\def\C{{\cal C}}
\def\D{{\bf D}}
\def\calS{{\cal S}}
\def\O{{\cal O}}
\def\P{{\cal P}}
\def\Q{{\cal Q}}
\def\u{{\cal U}}
\def\V{{\cal V}}
\def\F{{\cal F}}
\def\G{{\cal G}}
\def\H{{\cal H}}
\def\E{{\cal E}}
\def\m{{\bf m}}
\def\o{{\bf o}}
\def\LL{{\cal L}}
\def\MM{{\bf M}}
\def\DD{{\cal D}}
\def\RR{{\cal R}}
\def\NN{{\bf N}}
\def\EE{{\bf E}}
\def\FF{{\bf F}}
\def\GG{{\bf G}}
\def\T{{\cal T}}
\def\TT{{\bf T}}
\def\BB{{\bf B}}
\def\HH{{\bf H}}
\def\KP{{\bf KP}}
\def\K{{\bf K}}
\def\aa{{\bf a}}

\def\rank{{\rm Rank}}
\def\span{{\rm span}}
\def\dim{{\rm dim}}
\def\diam{{\rm diam}}
\def\diag{{\rm diag}}
\def\dist{{\rm dist}}
\def\trace{{\rm trace}}
\def\div{{\rm div}}
\def\spt{{\rm spt}}
\def\Lip{{\rm Lip}}
\def\Sp{{\bf S}}
\def\Z{{\bf Z}}

\def\Grass{{\rm Grass}}
\def\diam{{\rm diam}\,}
\def\rmB{{\rm B}}

\def\pa{\partial}
\def\um{^{-1}}
\def\pim{\pi^{-1}}
\def\pm{p^{-1}}
\def\cit{\cite{BouchitteBF2}}
\def\citd{\cite{GiaquintaM2}}
\def\ov{\overline}
\def\na{\nabla}
\def\vpi{\varpi}

\def\ti{\tilde}
\def\ess{{\rm ess}}
\def\ext{^{\rm ext}}
\def\u{{\bf u}}

\def\w{{\bf w}}
\def\e{{\bf e}}
\def\Om{\Omega}
\def\om{\omega}
\def\al{\aleph}

\def\ioe{\int_{\Omega_\veps}}
\def\ioesmbe{\int_{\Omega_\veps\setminus \C_\veps}}
\def\io{\int_{\Omega}}
\def\ise{\int_{D_\veps}}
\def\ibe{\int_{\C\veps}}
\def\iyke{\int_{Y^k_\veps}}
\def\dme{\;dm_\veps}
\def\Ome{{\Omega_\veps}}
\def\me{{m_{\veps}}}
\def\intb{{\int\!\!\!\!\!\!-}}

\def\QT{{\Omega^T}}
\def\De{{D_\veps}}
\def\Rex{{R_\veps}}
\def\Ri{{r_\veps}}

\title{\textbf{Homogenization of a diffusion process in a rarefied binary structure}}

\author{\textbf{Fadila~Bentalha $^*$, Isabelle~Gruais $^{**}$ and
Dan~Poli\v{s}evski $^{***}$} }
\date{}
\maketitle

{\bf Abstract.} We study the homogenization of a diffusion process which takes
place in a binary structure formed by an ambiental connected phase surrounding a
suspension of very small spheres distributed in an $\veps$-periodic network. The
asymptotic distribution of the concentration is determined for both phases, as
$\veps\to 0$, assuming that the suspension has   mass of
unity order and  vanishing volume.
Three
cases are distinguished according to the values of a certain limit capacity.
When it is positive and finite, the macroscopic system
involves a two-concentration system,  coupled through a  term  accounting for
the non local effects. In the   other two cases, where
the capacity is either infinite or going to zero, although the form of the system
is much simpler, some
peculiar effects still account for the presence of the suspension.



\section{Introduction}\lab{s:1}

\hspace{0.5cm}

Diffusion occurs naturally and is important in many industrial and geophysical
problems, particularly in oil recovery, earth pollution, phase transition,
chemical and nuclear processes. When one comes to
a rational study of  binary structures, a  crucial point lies in the interaction
between the microscopic and macroscopic levels and particularly the way the
former influences the latter. Once the distribution is assumed to be
$\veps$-periodic, this kind of study can be accomplished by the homogenization
theory.

The present study reveals the basic mechanism which governs
diffusion in both phases of such a binary structure, formed by an
ambiental connected phase surrounding a periodical suspension of
small particles. For simplicity, the particles are considered here
to be spheres of radius $r_\veps<<\veps$, that is
$\dis\lim_{\veps\to 0}\frac{r_\veps}{\veps} = 0$. We balance this
assumption, which obviously means that the suspension has vanishing
volume, by imposing the total mass of the suspension to be always of
unity order. This simplified structure permits the accurate
establishment of the macroscopic equations by means of a multiple
scale method of the homogenization theory adapted for fine-scale
substructures. It allows to have a general view on the specific
macroscopic effects which arise in every possible case. As we use
the non-dimensional framework, the discussion is made in fact with
respect to
only two parameters: $r_\veps$ and $b_\veps$, the latter standing   for the
ratio of suspension/ambiental phase diffusivities. As the
diffusivities of the two components can differ by orders of magnitude, the
interfacial  conditions play an important role.

It happens that the following cases have different treatments: $r_\veps<<\veps^3$,
$\veps^3 << r_\veps <<\veps$ and $r_\veps = \O(\veps^3)$.

To give a flavor of what may be considered as an appropriate choice
of the relative scales, we refer to the pioneering work
\cite{CioranescuM3} where the appearance of an extra term in the
limit procedure is responsible for a change in the nature of the
mathematical problem and is linked to  a critical size of the
inclusions. Later \cite{CasadoDiaz} showed how this could be
generalized to the $N$-dimensional case for non linear operators
satisfying classical properties of polynomial growth and coercivity.
Since then, the notion of non local effects has been developed in a
way that is closer to the present point of view in \cite{BellieudG},
\cite{BrianeT} and \cite{BGP}.

In dealing with our problem, the main difficulty was due to the
choice of test functions to be used in the associated  variational
formulation and which are classically some perturbation  of the
solution to the so-called cellular problem. Indeed, proceeding as
usual in homogenization theory, we use energy arguments based on a
priori estimates where direct limiting procedure apparently leads to
singular behavior. Non local effects appear when these singularities
can be overcome, which is usually achieved  by using adequate test
functions in the variational formulation. Since the fundamental work
\cite{CioranescuM3}, an important step was accomplished in this
direction in \cite{BellieudB}. A slightly different approach
\cite{BrianeT} uses Dirichlet forms involving non classical measures
in the spirit of \cite{Mosco}. However, the main drawback of this
method lies in its essential use of  the Maximum Principle, which
was avoided in \cite{BellieudG} for elastic fibers, and later in
\cite{BGP} where the case of spherical symmetry is solved. The
asymptotic behavior of highly heterogeneous media has also been
considered in the framework of homogenization when the coefficient
of one component is vanishing  and both components have volumes of
unity order: see the derivation of a double porosity model for a
single phase flow by \cite{ArbogastDouglasHornung} and the
application of two-scale convergence in order  to model diffusion
processes in \cite{Allaire4}.

The paper is organized as follows. Section~\ref{s:2} is devoted to
the main notations and to the description of the initial problem. We
set the functional framework (\ref{215}) where the existence and
uniqueness of the solution can be established: see  \cite{EneP1} and
\cite{Polisevski2} for similar problems. In Section~\ref{s:3}, we
introduce specific tools to handle the limiting  process. This is
based  on the use of the operators $G_r$ defined by  (\ref{gr})
which have a localizing effect: this observation motivates the
additional assumption (\ref{fwe2c}) on the external sources when the
radius of the particles is of critical order $\veps^3$ with $\veps $
denoting the period of the distribution. While passing to the limit,
the capacity number $\gamma_\veps$ defined by (\ref{gam})   appears
as the main criterium to describe the limit problem, the relative
parts played by the radius of the particles and by the period of the
network becoming explicit.

Section~\ref{s:4}, which is actually the most involving one,  deals
with the critical case when $\gamma_\veps$ has a  positive and
finite limit $\gamma$. In this part, where we assume also $b_\veps\to
+\infty$ the test functions are a convex combination of the elementary
solution of the Laplacian and its transformed by the operator
$G_\Ri$ defined by (\ref{gr}) with $r = \Ri$. This choice, which is
inspired  from \cite{BellieudB}, \cite{BellieudG} and \cite{BGP} and
has  to be compared with \cite{CioranescuM3}, allows to overcome the
singular behavior of the energy term when the period $\veps$ tends
to zero. We have to emphasize that this construction highly depends
on the geometry  of the problem, that is the spherical symmetry . To
our knowledge, the generalization to more intricate geometries
remains to be done. The resulting model
(\ref{pblim1})--(\ref{pblim2}), with the initial value defined after
$u_0$ in (\ref{cvu0}) and $v_0$ in (\ref{cvv0}), involves a pairing
$(u,v)$ which is coupled  through a linear operator acting on the
difference $u-v$ by the factor $4\pi\gamma$.

The case of the infinite capacity, where $\veps^3 << r_\veps <<
\veps$, is worked out in Section~\ref{s:5}. The proofs are only
sketched because  the arguments follow the same lines as in
Section~\ref{s:4}. Let us mention that the singular behavior of the
capacity  in this case, that is $\gamma_\veps\to+\infty$, forces $v$
to coincide  with  $u$. In other words, the infinite capacity
prevents the
splitting of the distribution, as it did in the critical case. Quite
interestingly, the initial value  of the global concentration is a convex
combination (\ref{inituinf}) of the initial conditions
$u_0$ and $v_0$; moreover, the mass density of the macroscopic diffusion equation
(\ref{pblim1inf}) takes both components into account, in accordance with the
intuition that the limiting process must lead to a binary mixture.

Finally, the case of vanishing capacity is handled in
Section~\ref{s:6}, that is when $r_\veps << \veps^3$. Here,  $v$
remains constant in time, obviously equal to the initial condition
$v_0$,  while $u$ satisfies the diffusion equation
(\ref{pblim1inf06})--(\ref{inituinf06})  with data independent of
the initial condition of the suspension. This can be seen as a proof
that when the radius of the particles is too small, then the
suspension does not present macroscopic effects, although a
corresponding residual  concentration, constant in time,  should be
considered.

\section{The diffusion problem}\lab{s:2}

We consider $\Om\subseteq\R^3$ a bounded Lipschitz domain occupied by a mixture
of two different materials, one of them forming the ambiental connected phase and
the other being concentrated in a periodical suspension of small spherical
particles.
Let us denote
\be
\lab{21}
 Y:= \left(- \frac{1}2,+\frac{1}2\right)^3.
\ee
\be
\lab{22}
 Y^k_\veps := \veps k + \veps Y,\quad k\in\Z^3.
\ee
\be
\lab{23}
 \Z_\veps := \{k\in\Z^3,\quad Y^k_\veps\subset \Om\},\quad
 \Om_{Y_\veps} := \cup_{k\in\Z_\veps}Y^k_\veps.
\ee

The suspension is defined by the following reunion
\be
\lab{24}
D_\veps := \cup_{k\in\Z_\veps}B(\veps k,r_\veps),
\ee
where $0 < r_\veps <<\veps$ and
$B(\veps k,r_\veps)$ is the ball of radius $r_\veps$ centered at
$\veps k$, $k\in\Z_\veps$. Obviously,
\be
\lab{de}
 \vert\De\vert\to 0\quad\mbox{as}\quad \veps\to 0.
\ee

The fluid domain is given by
\be
\lab{25}
 \Ome = \Om\setm D_\veps.
\ee

We also use the following notation
for the cylindrical time-domain:
\be
\lab{26}
 \Om^T := \Om\times ]0,T[;
\ee
similar definitions for $\Om^T_\veps$, $\Om_{Y_\veps}^T$ and $D^T_\veps$.

We  consider the problem which governs the diffusion process throughout our
binary mixture. Denoting by  $a_\veps>0$ and $b_\veps>0$  the relative
mass density and  diffusivity of the suspension, then, assuming without loss of
generality that
$\vert\Om\vert = 1$, its non-dimensional form is the following:

To find $u^\veps$
solution  of
\be
\lab{27}
 \rho^\veps\frac{\pa u^\veps}{\pa t} - \div{(k^\veps\na u^\veps)} =
f^\veps \quad\mbox{in}\quad
\Om^T
\ee
\be
\lab{28}
 [u^\veps]_\veps = 0 \quad\mbox{on}\quad\pa D^T_\veps
\ee
\be
\lab{29}
 [k^\veps\na u^\veps]_\veps n = 0 \quad\mbox{on}\quad\pa D^T_\veps
\ee
\be
\lab{210}
 u^\veps = 0\quad\mbox{on}\quad\pa\Om^T
\ee
\be
\lab{211}
 u^\veps(0) = u_0^\veps \quad\mbox{in}\quad\Om
\ee
where $[\cdot]_\veps$ is the jump across the interface $\pa D_\veps$,
$n$ is the normal on $\pa \De$ in the outward  direction, $f^\veps\in
L^2(0,T; H^{-1}(\Om))$, $u_0^\veps\in L^2(\Om)$ and

\be
\lab{212}
 \rho^\veps(x) = \left\{\ba{lll} 1 & \mbox{if} & x\in \Ome\\
a_\veps & \mbox{if} & x\in D_\veps
\ea\right.
\ee
\be
\lab{213}
 k^\veps(x) = \left\{\ba{lll} 1 & \mbox{if} & x\in \Ome\\
b_\veps & \mbox{if} & x\in D_\veps
\ea\right.
\ee

Let $H_\veps$ be the Hilbert space  $L^2(\Om)$ endowed with the scalar product
\be
\lab{214}
 (u,v)_{H_\veps} := (\rho^\veps u,v)_\Om
\ee
As $H^1_0(\Om)$ is dense in $H_\veps$ for any fixed $\veps>0$, we can set
\be
\lab{215}
 H^1_0(\Om)\subseteq H_\veps\simeq H_\veps'\subseteq H^{-1}(\Om)
\ee
with continuous embeddings.

Now, we can present the variational formulation of the problem
(\ref{27})-(\ref{211}).

To find $u^\veps\in L^2(0,T;H^1_0(\Om))\cap L^\infty(0,T;H_\veps)$
satisfying (in some sense) the initial condition (\ref{211}) and the following
equation
\be
\lab{216}
 \frac{d}{dt}(u^\veps,w)_{H_\veps} + (k_\veps\na u^\veps, \na w)_\Om =
\la f^\veps,w\ra \quad\mbox{in}\quad \DD'(0,T),\quad\forall w\in H^1_0(\Om)
\ee
where $\la\cdot,\cdot\ra$ denotes the duality product between $H^{-1}(\Om)$ and
$H^1_0(\Om)$.

\bthm
Under the above hypotheses and notations, problem (\ref{216}) has a unique
solution. Moreover, $\dis\frac{du^\veps}{dt}\in L^2(0,T;H^{-1}(\Om))$ and hence,
$u^\veps$ is equal almost everywhere to a function of $C^0([0,T]; H_\veps)$;
this is the sense of the initial condition (\ref{211}).
\ethm

In the following we consider that the density of the spherical particles is
much higher than that of the surrounding phase. The specific feature of our
mixture,  which describes the fact that although the
volume of the suspension is vanishing its  mass is of
unity order, is given by:

\be
\lab{ae}
 \lim_{\veps\to 0}a_\veps\vert\De\vert = a> 0
\ee

Regarding the relative diffusivity, we only assume:

\be
\lab{be}
  b_\veps \geq b_0>0,\quad\forall\veps>0.
\ee

As for the data, we assume that there exist  $f\in L^2(0,T; H^{-1}(\Om))$ and
$u_0\in L^2(\Om)$  such that
\be
\lab{fwe1}
 f^\veps\wto f\quad\mbox{in}\quad L^2(0,T; H^{-1}(\Om))
\ee
\be
\lab{cvu0}
 u_0^\veps\wto u_0\quad\mbox{in}\quad L^2(\Om)
\ee
Also, we assume that  there exist $C>0$ (independent of $\veps$) and $v_0\in
L^2(\Om)$ for which
\be
\lab{bornue0}
 \intb_\De \vert u_0^\veps\vert^2 dx\leq C
\ee
\be
\lab{cvv0}
 \frac{1}{\vert\De\vert}u_0^\veps\chi_{\De}\wto v_0\quad\mbox{in}\quad \DD'(\Om)
\ee
where, for any $D\subset\Om$, we denote
$$
 \intb_D \cdot dx = \frac{1}{\vert D\vert}\int_D\cdot dx.
$$

\brem
As $u^\veps_0$ satisfies (\ref{bornue0}) then  (\ref{cvv0}) holds
at least on some subsequence (see Lemma~A-2 \cite{BellieudB}).
\erem

\bprop
\lab{p:23}
We have
\be
\lab{ut}
 u^\veps\quad\mbox{is bounded in $L^\infty(0,T;L^2(\Om))\cap
L^2(0,T;H^1_0(\Om))$}.
\ee
Moreover, there exists $C>0$, independent of $\veps$, such that
\be
\lab{bornue}
 \intb_\De \vert u^\veps\vert^2 dx\leq C\quad\mbox{a.e. in $[0,T]$}
\ee
\be
\lab{c:qme}
b_\veps
\vert\na u^\veps\vert^2_{L^ 2(D^T_\veps)}\leq C.
\ee
\eprop
\bP
Substituting $w = u^\veps$ in
the variational problem (\ref{216}) and integrating over $(0,t)$ for any $t\in
]0,T[$, we get:
$$
 \frac{1}2\left( \vert u^\veps(t)\vert^2_\Ome + a_\veps
\vert u^\veps(t)\vert^2_\De\right)
 + b_\veps\int_0^t \vert\na u^\veps\vert^2_\De ds + \int_0^t\vert\na
u^\veps\vert^2_\Ome ds =
$$
$$
 = \int_0^t\la f^\veps(s), u^\veps(s)\ra ds + \frac{1}2\left(
\vert u_0^\veps\vert^2_\Ome + a_\veps\vert u_0^\veps\vert^2_\De\right).
$$

Notice that (\ref{cvu0}) and (\ref{bornue0}) yield:
$$
 \vert u_0^\veps\vert^2_\Ome + a_\veps\vert u_0^\veps\vert^2_\De\leq
\vert u_0^\veps\vert^2_\Om + a_\veps\vert\De\vert\intb_\De\vert u_0^\veps\vert^2
dx\leq C.
$$
Moreover:
$$
 \int_0^t\la f^\veps(s), u^\veps(s)\ra ds
\leq
\int_0^t\vert f^\veps\vert_{H^{-1}}\vert \na u^\veps\vert_\Om ds
$$
$$
 \leq
\int_0^t\vert f^\veps\vert_{H^{-1}}\vert \na u^\veps\vert_\Ome ds+
\int_0^t\vert f^\veps\vert_{H^{-1}}\vert \na u^\veps\vert_\De ds
$$
$$
 \leq \frac{1}2\int_0^T\!\!\!\vert f^\veps\vert_{H^{-1}}^2ds +
\frac{1}2\int_0^t\!\!\!\vert\na u^\veps\vert^2_\Ome ds+
\frac{1}{2b_\veps}\int_0^T\!\!\!\vert f^\veps\vert_{H^{-1}}^2 ds+
\frac{b_\veps}2\int_0^t\!\!\!\vert\na u^\veps\vert_\De^2 ds.
$$
There results:
$$
 \frac{1}2\left( \vert u^\veps(t)\vert^2_\Ome + a_\veps
\vert u^\veps(t)\vert^2_\De\right)
 +\frac{b_\veps}{2}\int_0^t \vert\na u^\veps\vert^2_\De ds +
\frac{1}2\int_0^t\vert\na u^\veps\vert^2_\Ome ds \leq
  C
$$
and the proof is completed.
\eP

\section{Specific tools}\lab{s:3}

First, we introduce
$$
 \RR_\veps = \{ R,\quad r_\veps << R <<\veps\}
$$
that is $R\in\RR_\veps$ iff
\be \lab{re} \quad \lim_{\veps\to 0}\frac{r_\veps}{R} =
\lim_{\veps\to 0} \frac{R}{\veps} = 0. \ee

We have to remark that $\RR_\veps$ is an infinite set, this property being
insured by the assumption $0 < r_\veps<<\veps$.

We denote the domain confined between the spheres of
radius
$a$ and
$b$ by
$$
 \C(a,b) :=
\{x\in\R^3,\; a < \vert x\vert< b\}
$$
and correspondingly
$$
 \C^k(a,b) := \veps  k + \C(a,b).
$$

For any $R_\veps\in\RR_\veps$, we use the following notations:
$$
 \C_\veps := \cup_{k\in\Z_\veps} \C^k(r_\veps,R_\veps),\quad
\C^T_\veps := \C_\veps\times ]0,T[
$$

\bd
For any $R_\veps\in\RR_\veps$, we define $w_{R_\veps}\in  H^1_0(\Om)$ by
\beq
\lab{defwe}
 w_{R_\veps}(x) &:=& \left\{\ba{c}
0\quad\mbox{in}\quad\Ome\setm \C_\veps,
\\
W_{R_\veps}( x - \veps k)\quad\mbox{in}\quad  \C^k_\veps,\quad
\forall k\in\Z_\veps,
\\
1\quad\mbox{in}\quad D_\veps.
\ea\right.
\eeq
where
\be
\lab{defbwe}
  W_{R_\veps}(y) =  \frac{r_\veps}{(R_\veps -
r_\veps)}\left(\frac{R_\veps}{\vert y\vert} - 1\right)\quad\mbox{for}\quad
y\in \C(r_\veps,R_\veps)
\ee

We have to remark here that $W_{R_\veps}\in H^1(C(r_\veps,R_\veps))$ and
satisfies the system
\beq
\lab{defbwe1}
 \Delta W_{R_\veps} & = & 0\quad\mbox{in}\quad \C(r_\veps,R_\veps)
\\
W_{R_\veps} & = & 1\quad\mbox{for}\quad  \vert y\vert = r_\veps
\\
\lab{defbwe2}
W_{R_\veps} & = & 0\quad\mbox{for}\quad  \vert y\vert = R_\veps
\eeq

\ed

From now on, we denote
\be
\lab{gam}
 \gamma_\veps := \frac{r_\veps}{\veps^3}.
\ee

\bprop
\lab{p:web}
For any $R_\veps\in\RR_\veps$, we have
\be
\lab{web}
 \vert\na w_{R_\veps}\vert_\Om\leq C\gamma_\veps^{1/2}
\ee
\be
\lab{we0}
w_{R_\veps}\to 0\quad\mbox{in}\quad L^2(\Om).
\ee
\eprop
\bP
First notice that
$$
 \vert w_{R_\veps}\vert_\Om = \vert w_{R_\veps}\vert_{\C_\veps\cup
D_\veps}\leq \vert \C_\veps\cup D_\veps\vert^{1/2} \leq
C\left(\frac{R_\veps}{\veps}\right)^{3/2}
$$
and $\lim_{\veps\to 0}\frac{R_\veps}{\veps} = 0$ by assumption (\ref{re}).

As for the rest, direct computation shows
\beqs
\vert\na w_{R_\veps}\vert_\Om^2
& = &
\sum_{k\in\Z_\veps}\int_{C^k_{r_\veps,R_\veps}}
\vert\na w_{R_\veps}\vert^2\;dx
\\
& = &
\sum_{k\in\Z_\veps}\int_0^{2\pi}\;d\Phi\int_0^{\pi}\sin\Theta\;d\Theta
\int_{r_\veps}^{R_\veps}\frac{dr}{r^2}\left(\frac{r_\veps
R_\veps}{R_\veps - r_\veps}\right)^2
\\
& \leq &
C\frac{\vert\Om\vert}{\veps^3}
\left(\frac{1}{r_\veps} - \frac{1}{R_\veps}\right)\left(\frac{r_\veps
R_\veps}{R_\veps - r_\veps}\right)^2
\leq C\frac{\gamma_\veps}{(1 - \frac{r_\veps}{R_\veps})}
\eeqs
and the proof is completed by (\ref{re}).
\eP

Lemmas~\ref{l:dr1r2} and \ref{l:ar} below are set without proof
since they are a three-dimensional  adaptation of
Lemmas~A.3 and A.4~\cite{BellieudB}.

\blem
\lab{l:dr1r2}
For every   $0 <r_1<r_2$ and $u\in H^1(C(r_1,r_2))$, the following estimate holds
true:
\be
\lab{c:dr1r2}
\vert\na u\vert^2_{C(r_1,r_2)}
\geq  \frac{4\pi r_1r_2}{r_2 - r_1}\left\vert
\intb_{\Sp_{r_2}}u\;d\sigma -
\intb_{\Sp_{r_1}}u\;d\sigma\right\vert^2,
\ee
where
$$
 \intb_{\Sp_r}\cdot\;d\sigma :=
\frac{1}{4\pi r^2}\int_{\Sp_r}\cdot\;d\sigma.
$$

\elem

\blem
\lab{l:ar}
There exists a positive constant $C>0$ such that: $\forall (R,\a)\in
\R^+\times (0,1)$, $\forall u\in H^1(B(0,R))$,
\be
\lab{ar}
 \int_{B(0,R)}\vert u- \intb_{\Sp_{\a R}}u\;d\sigma\vert^2\;dx\leq
C\frac{R^2}{\a}\vert\na u\vert^2_{B(0,R)}.
\ee
\elem

\bd

Consider the piecewise constant functions $G_r:\; L^2(0,T;H^1_0(\Om)) \to
L^2(\Om^T)$ defined  for any $r>0$ by
\be
\lab{gr}
G_r(\theta)(x,t)
=
\sum_{k\in\Z_\veps}\left(\intb_{\Sp^k_{r}}\theta(y,t)\;d\sigma_y
\right) 1_{Y^k_\veps}(x)
\ee
where we denote
\be
\lab{srk}
 S_{r}^k = \pa B(\veps k,r).
\ee

\ed

\blem
\lab{l:35}
If $R_\veps\in\RR_\veps$, then
for every $\theta\in  L^2(0,T;H^1_0(\Om))$ we have
\beq
\lab{I1}
 \vert\theta -
G_\Rex(\theta)\vert_{L^2(\Om^T_{Y_\veps})} & \leq &
C\left(\frac{\veps^3}{R_\veps}\right)^{1/2}\vert\na\theta\vert_{L^2(\QT)}
\\
\lab{I2}
\vert\theta - G_\Ri(\theta)\vert_{L^2(D^T_\veps)}
& \leq &
C r_\veps\vert\na\theta\vert_{L^2(D^T_\veps)}
\\
\lab{I3}
\vert G_\Rex(\theta) - G_\Ri(\theta)\vert_{L^2(\Om^T)}
&\leq &
C\left(\frac{\veps^3}{r_\veps}\right)^{1/2}\vert\na\theta\vert_{L^2(\C^T_\veps)}
\eeq
where $G_\Rex(\theta)$ and $G_\Ri(\theta)$ are defined following (\ref{gr}).

Moreover:
\be
\lab{mexb}
\!\!\!\!\!\!
 \vert G_\Rex(\theta)\vert^2_{L^2(\QT)} = \!\!
\int_0^T\!\!\!\intb_\De\!\!\!\vert G_\Rex(\theta)\vert^2 dxdt,\;
 \vert G_\Ri(\theta)\vert^2_{L^2(\QT)} \!=\!
\int_0^T\!\!\!\intb_\De\!\!\!\vert G_\Ri(\theta)\vert^2\;\!\!dxdt.
\ee

\elem
\bP
Notice that by definition:
$$
\sum_{k\in\Z_\veps}\int_0^T\!\!\!\iyke\vert \theta -
\intb_{\Sp^k_{R_\veps}}
\theta\;d\sigma\vert^2\;dxdt
\leq\sum_{k\in\Z_\veps}\int_0^T\!\!\!\int_{B(\veps
k,\frac{\veps\sqrt{3}}{2})}\vert
\theta -
\intb_{\Sp^k_{R_\veps}}
\theta\;d\sigma\vert^2\;dxdt
$$
where we have used that
$$
 Y^k_\veps\subset B(\veps k,\frac{\veps\sqrt{3}}{2})
$$
for every $k\in\Z_\veps$. We use Lemma~\ref{l:ar} with
$$
 R = \frac{\veps\sqrt{3}}{2},\quad \a = \frac{2R_\veps}{\veps\sqrt{3}}
$$
to deduce that
$$
\int_{\Om^T_{Y_\veps}}\vert\theta - G_\Rex(\theta)\vert^2\;dxdt
\leq
C\left(\frac{\veps\sqrt{3}}{2}\right)^2\frac{\veps\sqrt{3}}{2R_\veps}
\sum_{k\in\Z_\veps}\int_0^T\!\!\!
\int_{B(\veps k,\frac{\veps\sqrt{3}}{2})}\vert\na\theta\vert^2\;dxdt
$$
$$
\leq
C\frac{\veps^3}{R_\veps}\sum_{k\in\Z_\veps}\int_0^T\!\!\!
\int_{B(\veps k,\frac{\veps\sqrt{3}}{2})}\vert\na\theta\vert^2\;dxdt
\leq
C\frac{\veps^3}{R_\veps}\int_\QT\vert\na\theta\vert^2\;dxdt
$$
which shows (\ref{I1}).

To establish (\ref{I2}), we recall the definition:
$$
\int_{D^T_\veps}\vert\theta - G_\Ri(\theta)\vert^2\;dxdt
=
\sum_{k\in\Z_\veps}\int_0^T\!\!\!\int_{B(\veps k,r_\veps)}\vert \theta -
\intb_{\Sp^k_{r_\veps}}\theta\;d\sigma\vert^2\;dxdt
$$
Applying Lemma~\ref{l:ar} with $R = r_\veps$ and $\a = 1$,
we get the result
\beqs
\int_{D^T_\veps}\vert\theta - G_\Ri(\theta)\vert^2\;dxdt
&\!\!\! \leq\!\!\! &
Cr_\veps^2\sum_{k\in\Z_\veps}\int_0^T\!\!\!\int_{B(\veps k,r_\veps)}
\!\!\!\vert\na\theta\vert^2\;dxdt
\leq C r_\veps^2\int_{D^T_\veps}\vert\na\theta\vert^2\,dxdt.
\eeqs

We come to (\ref{I3}). Indeed, applying Lemma~\ref{l:dr1r2} and
(\ref{re}):
$$
\int_\QT\vert G_\Rex(\theta) - G_\Ri(\theta)\vert^2\;dxdt
=
\sum_{k\in\Z_\veps}\int_0^T\!\!\!\iyke\vert\intb_{\Sp^k_{R_\veps}}
\theta\;d\sigma
- \intb_{\Sp^k_{r_\veps}}\theta\;d\sigma\vert^2\;dydt
$$
$$
\leq
\sum_{k\in\Z_\veps}\iyke\frac{(R_\veps - r_\veps)}{4\pi R_\veps r_\veps}
\;dy\int_0^T\!\!\!\int_{C^k_{r_\veps,R_\veps}}\!\!\!
\vert\na\theta\vert^2\;dxdt
=
\frac{(R_\veps - r_\veps)}{4\pi r_\veps R_\veps}
\sum_{k\in\Z_\veps}\veps^3\int_0^T\!\!\!
\int_{C^k_{r_\veps,R_\veps}}\!\!\!\vert\na\theta\vert^2\,dxdt
$$
$$
=
C\veps^3\frac{(R_\veps - r_\veps)}{4\pi r_\veps R_\veps}\int_{\C^T_\veps}
\vert\na\theta\vert^2\;dxdt
\leq
C\frac{\veps^3}{r_\veps}\int_{\C^T_\veps}\vert\na\theta\vert^2\;dxdt.
$$
Finally, a direct computation yields (\ref{mexb}).
\eP

\bprop
\lab{p:naom}
If $R_\veps\in\RR_\veps$, then
for any $\theta\in L^2(0,T;H^1_0(\Om))$ there holds true:
$$
   \int_0^T\!\!\!\intb_\De\vert\theta\vert^2\;dxdt\leq
C\max{(1,\frac{\veps^3}{r_\veps})}
\vert\na\theta\vert^2_{L^2(\QT)}.
$$
\eprop
\bP
We have:
$$
\int_0^T\!\!\!\intb_\De\vert\theta\vert^2\;dxdt
\leq
2\int_0^T\!\!\!\intb_\De\vert\theta - G_\Ri(\theta)\vert^2\;dxdt +
2 \int_0^T\!\!\!\intb_\De\vert G_\Ri(\theta)\vert^2\;dxdt
$$
$$
=
2\int_0^T\!\!\!\intb_\De\vert\theta - G_\Ri(\theta)\vert^2\;dxdt +
2 \int_\QT\vert G_\Ri(\theta)\vert^2\;dxdt
$$
$$
\leq
Cr^2_\veps\int_0^T\!\!\!\intb_\De\vert\na\theta\vert^2\;dxdt +
4 \int_\QT\vert G_\Ri(\theta) - G_\Rex(\theta)\vert^2\;dxdt +
$$
$$
+
8 \int_\QT\vert G_\Rex(\theta) - \theta\vert^2\;dxdt +
8 \int_\QT\vert\theta\vert^2\;dxdt
$$
$$
\leq
Cr^2_\veps\int_0^T\!\!\!\intb_\De
\vert\na\theta\vert^2\;dxdt
+  C \frac{\veps^3}{r_\veps}
\int_{\C^T_\veps}\vert\na\theta\vert^2\;dxdt +
$$
$$
+
C\frac{\veps^3}{R_\veps}\int_\QT\vert\na\theta\vert^2\;dxdt
+
C\int_\QT\vert\na\theta\vert^2\;dxdt
$$
$$
\leq
C\left(\frac{\veps^3}{r_\veps} + \frac{\veps^3}{R_\veps} + 1\right)
\int_\QT\vert\na\theta\vert^2\;dxdt
\leq C\max{(1,\frac{\veps^3}{r_\veps})}\int_\QT\vert\na\theta\vert^2\;dxdt
$$
\eP

\brem
\lab{r:37}
Using the Mean Value Theorem, we easily find that
$$
 \vert G_\Ri(\varphi) - \varphi\vert_{L^\infty(\C_\veps\cup\De)}\leq
2R_\veps\vert\na\varphi\vert_{L^\infty(\Om)},\quad\forall\varphi\in\DD(\Om).
$$
\erem

\bd
Let $M_{D_\veps}:\; L^2(0,T;C_c(\Om))\to L^2(\QT)$ be defined by
$$
 M_\De(\varphi)(x,t) :=
\sum_{k\in\Z_\veps}\left(\intb_{Y^k_\veps}\varphi(y,t)\;dy
\right)\;1_{B(\veps k,r_\veps)}(x).
$$
\ed

\blem
\lab{l:ffe}
For any $\varphi\in L^2(0,T;C_c(\Om))$, we have:
$$
 \lim_{\veps\to 0}\int_0^T\!\!\!\intb_\De\vert\varphi - M_\De(\varphi)\vert^2dxdt
= 0.
$$
\elem
\bP
Notice that
$$
 \intb_\De\vert\varphi - M_{D_\veps}(\varphi)\vert^2
dx  =
\frac{1}{\vert \De\vert}\sum_{k\in\Z_\veps}
\int_{B(\veps k,r_\veps)}\vert\varphi - \intb_{Y^k_\veps}\varphi\;dy
\vert^2\;dx.
$$
As
$
{\rm card}(\Z_\veps) \simeq
\dis\frac{\vert\Om\vert}{\veps^3},\quad\mbox{then}\quad
\dis\vert B(0,r_\veps)\vert \frac{{\rm card}(\Z_\veps)}{\vert\De\vert}\to
\vert\Om\vert = 1
$
and by the uniform continuity of $\varphi$ on $\Om$ it follows the convergence to
$0$ a.e. on $[0,T]$. Lebesgue's dominated convergence theorem achieves the result.
\eP

\section{Homogenization of the case $r_\veps = \O(\veps^3)$}\lab{s:4}

The present critical radius case is described by
\be
\lab{cap}
\lim_{\veps\to 0}\gamma_\veps =  \gamma\in ]0,+\infty[.
\quad
\ee
Its homogenization process is the most involving one. That is why we start
the homogenization study of our problem with this case, under the
condition
\be
\lab{binf}
 \lim_{\veps\to 0}b_\veps = +\infty
\ee

We also assume  that $f^\veps$ has the following additional  property:
\be
\lab{fwe2c}
\left\{
\ba{l}
 \exists R_\veps\in\RR_\veps
\quad\mbox{ and}\quad  g\in L^2(0,T;H^{-1}(\Om))\quad\mbox{for which}\quad
\\
\\
 \la f^\veps,w_{R_\veps}\varphi\ra\to \la g,\varphi\ra
\quad\mbox{in}\quad\DD'(0,T),\quad\forall
\varphi\in\DD(\Om)
\ea\right.
\ee
(see \cite{BGP}
for a certain type of functions  $f^\veps$ which satisfy (\ref{fwe2c})).

\brem
Notice that due to (\ref{cap}), in this case Proposition~\ref{p:naom}  reads
\be
\lab{tme1}
\forall\varphi\in L^2(0,T;H^1_0(\Om)),\quad
 \int_0^T\!\!\!\intb_\De\vert\varphi\vert^2\;dxdt\leq
C\vert\na\varphi\vert^2_{L^2(\QT)}.
\ee
\erem

A preliminary result is the following:

\bprop
\lab{p:uttau}
There exist $ u\in L^\infty(0,T;L^2(\Om))\cap L^2(0,T;H^1_0(\Om))$ and $v\in
L^2(\QT)$ such that, on some subsequence,
\be
\lab{cvstar}
u^\veps  \buildrel{\star}\over\wto  u\quad\mbox{in}\quad L^\infty(0,T;L^2(\Om))
\ee
\be
\lab{cvl2}
u^\veps  \wto  u\quad\mbox{in}\quad L^2(0,T;H^1_0(\Om))
\ee
\be
\lab{cvgre}
G_\Rex(u^\veps)\to u\quad\mbox{in}\quad L^2(\QT)
\ee
\be
\lab{cvgri}
G_\Ri(u^\veps) \wto v\quad\mbox{in}\quad
L^2(\QT)
\ee
Moreover, we have
\be\lab{cvugre}
\lim_{\veps\to 0}
\int_0^T\!\!\!\intb_\De \vert u^\veps - G_{r_\veps}(u^\veps)\vert^2 dxdt = 0
\ee
\eprop
\bP
From (\ref{ut}), we get, on some subsequence, the  convergences (\ref{cvstar})
and (\ref{cvl2}).
Moreover, we have:
\be
\lab{ugre1}
 \vert u-G_\Rex(u^\veps)\vert^2_\QT = \vert u\vert^2_{\QT\setm\Om^T_{Y_\veps}}  +
\vert u-G_\Rex(u^\veps)\vert^2_{\Om^T_{Y_\veps}}
\ee
where:
\be
\lab{ugre}
 \vert u-G_\Rex(u^\veps)\vert_{\Om^T_{Y_\veps}} \leq
\vert u-u^\veps\vert_{\Om^T_{Y_\veps}}  + \vert
u^\veps-G_\Rex(u^\veps)\vert_{\Om^T_{Y_\veps}}
\ee
$$
 \leq \vert u-u^\veps\vert_{\QT} + \vert
u^\veps-G_\Rex(u^\veps)\vert_{\Om^T_{Y_\veps}}
$$
and (\ref{I1}) yields:
$$
 \vert u^\veps-G_\Rex(u^\veps)\vert^2_{\Om^T_{Y_\veps}} \leq
C\frac{\veps^3}{R_\veps}
\vert\na u^\veps\vert^2_\QT =
C\frac{\veps^3}{r_\veps}\frac{r_\veps}{R_\veps}
\vert\na u^\veps\vert^2_\QT\leq C\frac{r_\veps}{R_\veps}
$$
and thus:
$$
 \lim_{\veps\to 0}\vert u^\veps-G_\Rex(u^\veps)\vert^2_{\Om^T_{Y_\veps}}  = 0.
$$
As (\ref{cvl2}) implies that
\be
\lab{strcv}
 u^\veps\to u\quad\mbox{in}\quad L^2(\QT)
\ee
the right-hand side of (\ref{ugre}) tends to zero as $\veps\to 0$,  that is:
$$
 \lim_{\veps\to 0}\vert u-G_\Rex(u^\veps)\vert_{\Om^T_{Y_\veps}} = 0.
$$
After substitution into the right-hand side of (\ref{ugre1}), and taking into
account that
$$
 \lim_{\veps\to 0} \vert\QT\setm\Om^T_{Y_\veps}\vert = 0,
$$
we obtain (\ref{cvgre}), that is,
\be
\lab{ttild}
G_\Rex(u^\veps)
 \to u\quad\mbox{in}\quad L^2(\QT).
\ee

In order to prove (\ref{cvgri}), we see that
\be
\lab{ttc}
\ba{l}
 \vert G_\Ri(u^\veps)\vert_{L^2(\QT)}
\leq \dis
\vert G_\Ri(u^\veps) - G_\Rex(u^\veps)\vert_{L^2(\QT)} +
\vert G_\Rex(u^\veps)\vert_{L^2(\QT)}
\\ \\ \dis
\leq
\frac{^C}{\gamma^{1/2}_\veps}
\vert\na u^\veps\vert_{L^2(\QT)} + C
\leq C.
\ea
\ee

Moreover, recall that from (\ref{I2}) we have, taking into account
(\ref{c:qme}):
\be
\lab{qtt}
 \int_0^T\!\!\!\intb_\De\vert u^\veps - G_\Ri(u^\veps)\vert^2\;dx dt\leq C
r^2_\veps
\int_0^T\!\!\!\intb_\De\vert\na u^\veps\vert^2\;dx dt\leq \frac{C}{\gamma_\veps
b_\veps}\to 0
\ee
and the proof is completed.
\eP

\bprop
\lab{p:cvmue}
For any $\varphi\in L^2(0,T; C_c(\Om))$, we have:
\be
\lab{cvm}
\lim_{\veps\to 0}\int_0^T\!\!\!\intb_\De u^\veps\varphi dxdt=
\int_\QT v\varphi dxdt.
\ee
\eprop
\bP
We have:
\be
\lab{tt0}
\ba{l}\dis
 \int_0^T\!\!\!\intb_\De
u^\veps\varphi dxdt = \int_0^T\!\!\!\intb_\De
(u^\veps -  G_\Ri(u^\veps))\varphi dxdt +
\\
\\ \dis
+\int_0^T\!\!\!\intb_\De
G_\Ri(u^\veps) (\varphi - M_\De(\varphi)) dxdt
+ \dis \int_0^T\!\!\!\intb_\De G_\Ri(u^\veps) M_\De(\varphi) dxdt
\ea
\ee
The first right-hand term tends to zero thanks to (\ref{cvugre}) in
Proposition~\ref{p:uttau}. The second one tends also to zero thanks to
Lemma~\ref{l:ffe}.
The last term is handled as follows:
$$
 \int_0^T\!\!\!\intb_\De G_\Ri(u^\veps) M_\De(\varphi) dxdt=
\l_\veps
\sum_{k\in\Z_\veps}\int_0^T\!\!\!\int_{Y^k_\veps}
\left(\intb_{\Sp^k_{\veps}}u^\veps d\sigma\right) \varphi dx dt
 = \l_\veps\int_\QT\varphi G_{r_\veps}(u^\veps) dxdt
$$
where
$$
 \l_\veps :=\frac{\vert B(0,r_\veps)\vert}{\veps^3\vert\De\vert}\to
1\quad\mbox{as}\quad \vert\Om\vert = 1.
$$
The proof is completed by (\ref{cvgri}).
\eP

\bprop
\lab{c:44}
For any $\varphi\in L^2(0,T; H^1_0(\Om))$, we have
\be
\lab{ccvm}
 \int_0^T\!\!\!\intb_\De u^\veps\varphi dxdt\to \int_\QT v\varphi dxdt.
\ee
\eprop
\bP
In the light of proposition~\ref{p:cvmue}, we have to prove that the left-hand
side term is continuous in the corresponding norm. This can be obtained as
follows:
$$
 \left\vert\int_0^T\!\!\!\intb_\De u^\veps\varphi dxdt\right\vert
 \leq \left(\int_0^T\!\!\! \intb_\De\vert u^\veps\vert^2dxdt\right)^{1/2}
\left( \int_0^T\!\!\!\intb_\De\vert \varphi\vert^2dxdt\right)^{1/2} \leq
$$
$$
\leq  C\vert\varphi\vert^2_{L^2(0,T;H^1_0(\Om))},
$$
where we used (\ref{bornue}) and  (\ref{tme1}).
\eP

\bprop
\lab{p:fie}
Let for any $R_\veps\in\RR_\veps$ and $\varphi,\psi\in\DD(\Om)$
\be
\lab{deffie}
\Phi^\veps = (1-w_{R_\veps})\varphi + w_{R_\veps} G_\Ri(\psi)
\ee
Then, for any $\eta\in\DD([0,T[)$, we have
\be
\lab{cvfie}
 \lim_{\veps\to 0}\vert\Phi^\veps - \varphi\vert_\Om = 0
\ee
\be
\lab{lvar}
 \lim_{\veps\to 0} \int_\QT \rho^\veps u^\veps\Phi^\veps\eta'(t) dxdt=
\int_\QT u\varphi\eta'(t) dxdt + a\int_\QT v\psi\eta'(t) dxdt. \ee
\eprop \bP The property (\ref{cvfie}) is an immediate consequence of
(\ref{we0}) and of the uniform boundness of $G_\Ri(\psi)$ in
$L^\infty(\Om)$.

For the second property, let us notice that
$$
 \int_\QT \rho^\veps u^\veps\Phi^\veps\eta'(t) dxdt =
\int_0^T\!\!\!\int_\Om\chi_{\Ome} u^\veps\Phi^\veps(x)\eta'(t)dxdt
$$
$$ +
a_\veps\int_0^T\int_\De u^\veps G_\Ri(\psi)\eta'(t) dxdt.
$$
As we obviously have
$$
 \lim_{\veps\to 0}\int_0^T\!\!\!\int_\Om\chi_{\Ome}
u^\veps\Phi^\veps(x)\eta'(t)dxdt =
\int_\QT u\varphi\eta'(t) dxdt,
$$
it remains to study
$$
 a_\veps\int_0^T\int_\De u^\veps G_\Ri(\psi)\eta'(t) dxdt=
a_\veps\vert\De\vert\int_0^T\!\!\!\intb_\De u^\veps G_\Ri(\psi)\eta'(t) dxdt.
$$

Using (\ref{cvm}) and the uniform continuity of $\psi$,  we get
$$
 \lim_{\veps\to 0}
a_\veps\int_0^T\int_\De u^\veps G_\Ri(\psi)\eta'(t) dxdt=
a \int_\QT v\psi\eta' dxdt.
$$
\eP

\bprop
\lab{p:vare}
If  $\Phi^\veps$ is defined like in Proposition~\ref{p:fie}, then  we have
\be
\lab{gugv}
 \lim_{\veps\to 0}
\int_0^T\!\!\!\int_\Om \na u^\veps\cdot
\na\Phi^\veps \eta(t)\;dxdt
=\int_\QT \!\!\!\na u\cdot\na\varphi \eta(t)\;dxdt +
 4\pi\gamma\!\!
\int_\QT\!\!\!(v - u)(\psi -\varphi) \eta(t)\;dxdt
\ee
\eprop
\bP
First consider%
$$
 \int_0^T\!\!\!\int_\Ome \na u^\veps\cdot
\na\Phi^\veps\;dxdt
$$
which reduces to
$$
  \int_0^T\!\!\!\int_{\Ome\setm\C_\veps}  \na u^\veps\cdot
\na\varphi\eta\;dxdt +
\int_0^T\!\!\!\int_{\C_\veps}  \na u^\veps\cdot\na\Phi^\veps\eta dxdt.
$$
Lebesgue's dominated convergence theorem yields
$
 \na\varphi 1_{\Ome\setm \C_\veps}\to \na\varphi
$
in $L^2(\Om)$. Thus, taking (\ref{cvl2}) into account
$$
 \int_0^T\io\na u^\veps\cdot\na\varphi \,\eta \chi_{\Ome\setm
\C_\veps}\;dxdt\to
\int_\QT\na u\cdot\na\varphi\,\eta\;dxdt.
$$
Now, we come to the remaining part, namely
\be
\lab{be0}
\ba{lcl}
&&
 \dis\int_0^T\ibe \na u^\veps\cdot
\dis\na\Phi^\veps\eta(t)\;dxdt
=
 \int_0^T\ibe(1 - w_\Rex) \na u^\veps\cdot
\na\varphi\eta \;dxdt
\\
\!\!\!\!\!
&&\!\!\!\!\!+
\dis \int_0^T\!\!\!\ibe \!\!\na u^\veps\!\!\cdot
\na w_{R_\veps} (G_\Ri(\psi) - \varphi)dxdt
\\
& := & I_1 + I_2
\ea
\ee
In the first integral, as $\chi_{\C_\veps}\na\varphi\to 0$ in $L^2(\QT)$,
$\na u^\veps\wto\na u$ in $L^2(\QT)$ and $(1 - w_\Rex)$ is obviously
bounded, we easily find  that
$I_1$ tends to zero.

In order to study   $I_2$, let us notice that
\be
\lab{be5}
\ba{lcl}
I_2 &=& \!\!\!
\dis \int_0^T\!\!\!\ibe \!\!\!\! \na u^\veps\cdot
\na w_{R_\veps} (G_\Ri(\varphi) - \varphi)\,\eta \;dxdt+
\\
&&
\quad\dis +\int_0^T\!\!\!\ibe \!\!\!\!\na u^\veps\cdot
\na w_{R_\veps} (G_\Ri(\psi) - G_\Ri(\varphi))\,\eta \;dxdt
\ea
\ee
The first term in the
right-hand side of (\ref{be5}) may be estimated by
\be
\lab{be6}
\!\!\!
 \vert \int_0^T\!\!\!\ibe\!\!\! \na u^\veps\cdot
\na w_{R_\veps} (\varphi - G_\Ri(\varphi))\,\eta\;dx dt\vert\leq
\vert\na u^\veps\vert_\QT\vert\na w_{R_\veps}\eta\vert_\QT\vert\varphi -
G_\Ri(\varphi)\vert_{L^\infty(\C_\veps)} .
\ee
As $(w_{R_\veps})$ is bounded in $H^1(\Om)$
({\it see} Proposition~\ref{p:web}),
the right hand side
of (\ref{be6}) tends to zero by Remark~\ref{r:37}.

Going back to the second term in
the right hand side of (\ref{be5}), we may write
\beqs
&&
\int_0^T\!\!\!\ibe \na u^\veps\cdot
\na w_{R_\veps} (G_\Ri(\psi) - G_\Ri(\varphi))\eta(t)\;dxdt
\\
& = & \!\!\!\!
\sum_{k\in\Z_\veps}
\left(\!\intb_{\Sp^k_{r_\veps}}\!\!\psi\;d\sigma-
\intb_{\Sp^k_{r_\veps}}\!\!\varphi\;d\sigma\!\!\right)\int_0^{2\pi}\!\!
d\Phi\!\!\int_{0}^{\pi}\!\!\!\sin\Theta\;d\Theta
\int_{r_\veps}^{R_\veps}\left(\left.\!\int_0^T\frac{\pa u^\veps}{\pa
r} \right\vert_{\C^k(r_\veps,R_\veps)} \eta(t)dt\right) \!
\frac{d W_{R_\veps}}{dr} r^2\;dr\!\!
\\
& = &\!\!\!
\frac{r_\veps R_\veps}{(R_\veps - r_\veps)}
\sum_{k\in\Z_\veps}
\left(\intb_{\Sp^k_{r_\veps}}\psi\;d\sigma-
\intb_{\Sp^k_{r_\veps}}\varphi\;d\sigma\right)\int_{\Sp_1}
\int_0^T\!\!\!( u^\veps \vert_{\vert x - \veps
k\vert = r_\veps} -  u^\veps \vert_{\vert x - \veps k\vert =
R_\veps})\eta(t)dt
d\sigma_1
\\
& = &
\frac{4\pi r_\veps R_\veps}{\veps^3(R_\veps - r_\veps)}
\int_\QT (G_\Ri(u^\veps) - G_\Rex(u^\veps))(G_\Ri(\psi)
-G_\Ri(\varphi))\eta(t)\;dxdt
\eeqs
from which we infer that $I_2$ is converging to
$$
4\pi\gamma
\int_\QT(v - u)(\psi -\varphi)\eta(t)\;dxdt
$$
and the proof is completed.
\eP

We are in the position to state our  main result:

\bthm
\lab{c:lim}
The limits  $u\in L^\infty(0,T;L^2(\Om))\cap L^2(0,T;H^1_0(\Om))$ and $v\in
L^2(\QT)$  of (\ref{cvstar})--(\ref{cvgri})  verify (in a weak sense) the
following problem:
\beq
\lab{pblim1}
\frac{\pa u}{\pa t}
-\Delta u + 4\pi\gamma (u - v) & = &
(f-g)\quad\mbox{in}\quad\QT,
\\
a\frac{\pa v}{\pa t} + 4\pi\gamma (v - u) & = &
g\quad\mbox{in}\quad\QT,
\\
\lab{initu}
u(0) &= & u_0\quad\mbox{in}\quad\Om
\\
\lab{pblim2}
v(0) & = & v_0\quad\mbox{in}\quad\Om
\eeq
Moreover, there holds $u\in C^0([0,T]; L^2(\Om))$ and $v\in C^0([0,T];
H^{-1}(\Om))$; these are the senses of (\ref{initu}) and (\ref{pblim2}).
\ethm
\brem
As the problem (\ref{pblim1})--(\ref{pblim2}) has a unique weak solution,
the convergences in Proposition~\ref{p:uttau} hold on the whole sequence.
\erem

\bPof{Theorem \ref{c:lim}}
We set in (\ref{216}) $w = \Phi^\veps$ where $\Phi^\veps$ is defined like in
lemma~\ref{p:fie}. Then, by multiplying (\ref{216}) by $\eta\in
\DD([0,T[)$  and integrating it over $[0,T]$ we get
\be
\lab{var}
\!\!\!
 \dis -\int_\QT \!\!\rho^\veps u^\veps\Phi^\veps\eta'dxdt +
\int_\QT\!\! k^\veps\na u^\veps (\na \Phi^\veps)\eta dxdt =
\int_0^T\!\!\!\la f^\veps,\Phi^\veps\ra\eta dt + \int_\Om \!\rho^\veps u_0^\veps
\Phi^\veps\eta(0) dx.
\ee
Then, the left-hand side tends to
\be
\lab{c:vare}
\ba{c}
\dis
- \int_\QT\!\!\! u\varphi\eta' dxdt - a\int_\QT\!\!\! v\varphi\eta' dxdt +
 \int_\QT\!\!\! \na u\cdot\na\varphi\,\eta\;dxdt +
\\
\\
+\dis
 4\pi\gamma
\int_\QT\!\!\! (v - u)(\psi -\varphi)\,\eta\;dxdt.
\ea
\ee

This is a direct consequence of Proposition~\ref{p:vare} together with
the remark that
$$
 \int_0^T\!\!\!\int_{D_\veps}\na u^\veps\na\Phi^\veps\;dxdt = 0
$$
since $\Phi^\veps$ is constant on every $B(\veps k,r_\veps)$,
$k\in\Z_\veps$.

As for the right-hand side, we have
$$
 \int_0^T\la f^\veps,\Phi^\veps\ra\eta dt =
\int_0^T\la f^\veps,(1-w_{R_\veps})\varphi\ra\eta dt +  \int_0^T\la
f^\veps,w_\Rex G_\Ri (\psi)\ra\eta dt
$$
and, with hypothesis (\ref{fwe2c}),
$$
 \int_0^T\la f^\veps,(1-w_{R_\veps})\varphi\ra\eta dt \to
\int_0^T\la f-g,\varphi\ra\eta dt.
$$
Moreover,
$$
 \int_0^T\la
f^\veps,w_\Rex G_\Ri (\psi)\ra\eta dt =
\int_0^T\la
f^\veps,w_\Rex (G_\Ri (\psi)-\psi)\ra\eta dt +\int_0^T\la
f^\veps,w_\Rex \psi\ra\eta dt
$$
with
$$
 \left\vert
\int_0^T\la
f^\veps,w_\Rex (G_\Ri (\psi)-\psi)\ra\eta dt \right\vert\leq
\int_0^T\vert f^\veps\vert_{H^{-1}} \vert
w_\Rex(G_\Ri(\psi)-\psi)\vert_{H^1_0(\Om)}.
$$
As we have
$$
 \vert w_\Rex(G_\Ri(\psi)-\psi)\vert_{H^1_0(\Om)} =
\vert\na (w_\Rex(G_\Ri(\psi)-\psi))\vert_\Om
$$
$$
\leq
\vert\na w_\Rex\vert_{\C_\veps} \vert G_\Ri(\psi)-\psi\vert_{L^\infty(\C_\veps)}
+
\vert\na\psi\vert_{\C_\veps\cup\De}
$$
Remark~\ref{r:37} and (\ref{web}) obviously yield
$$
 \lim_{\veps\to 0}\vert w_\Rex(G_\Ri(\psi)-\psi)\vert_{H^1_0(\Om)} = 0.
$$
The assumption (\ref{fwe1}) on $f^\veps$ implies that $\vert
f^\veps\vert_{H^{-1}}\leq C$ and thus
$$
 \lim_{\veps\to 0}\int_0^T\la
f^\veps,w_\Rex (G_\Ri (\psi)-\psi)\ra\eta dt = 0.
$$
We conclude thanks to hypothesis (\ref{fwe2c}) that
$$
\lim_{\veps\to 0} \int_0^T\la
f^\veps,w_\Rex \psi\ra\eta dt = \int_0^T\la
g,\, \psi\ra\eta dt.
$$
Finally:
$$
 \lim_{\veps\to 0} \int_0^T\la f^\veps,\Phi^\veps\ra\eta dt =
\int_0^T\la
f-g,\, \varphi\ra\eta dt+ \int_0^T\la
g,\, \psi\ra\eta dt.
$$

We get
$$
  \int_\Om \!\rho^\veps u_0^\veps
\Phi^\veps\eta(0) dx =
\int_\Ome u_0^\veps\Phi^\veps\eta(0) dx + a_\veps\int_\De u_0^\veps
G_\Ri(\psi)\eta(0) dx.
$$
Using the hypotheses (\ref{cvu0})--(\ref{cvv0}) on
$u_0^\veps$, we  pass to the limit and with the same arguments as above we obtain
$$
 \lim_{\veps\to 0}\int_\Om \!\rho^\veps u_0^\veps
\Phi^\veps\eta(0) dx =
\eta(0)
\int_\Om (u_0\varphi + av_0\psi) dx
$$
which achieves the proof.
\ePof

\section{Homogenization in the case $\veps^3 << r_\veps <<\veps$}\lab{s:5}

In this section, we fix some $R_\veps\in\RR_\veps$.

\brem
\lab{r:51}
Notice that in this case Proposition~\ref{p:naom} also reads
\be
\lab{5:tme1}
 \int_0^T\!\!\!\intb_\De\vert\varphi\vert^2\;dxdt\leq
C\vert\na\varphi\vert^2_{L^2(\QT)},\quad
\forall\varphi\in L^2(0,T;H^1_0(\Om)).
\ee
\erem

In the present case, Proposition~\ref{p:23} and Lemma~\ref{l:35} imply in a
straightforward manner the result corresponding to Proposition~\ref{p:uttau}.

\bprop
\lab{p:uttauinf}
There exists $ u\in L^\infty(0,T;L^2(\Om))\cap L^2(0,T;H^1_0(\Om))$  such that,
on some subsequence,
\be
\lab{cvstari}
u^\veps  \buildrel{\star}\over\wto  u\quad\mbox{in}\quad L^\infty(0,T;L^2(\Om))
\ee
\be
\lab{cvl2i}
u^\veps  \wto  u\quad\mbox{in}\quad L^2(0,T;H^1_0(\Om))
\ee
\be
\lab{cvgrei}
G_\Rex(u^\veps)\to u\quad\mbox{in}\quad L^2(\QT)
\ee
\be
\lab{cvgrii}
G_\Ri(u^\veps) \to u\quad\mbox{in}\quad
L^2(\QT)
\ee
Moreover, we have
\be\lab{cvugrei}
\lim_{\veps\to 0}
\int_0^T\!\!\!\intb_\De \vert u^\veps - G_{r_\veps}(u^\veps)\vert^2 dxdt = 0
\ee
\eprop

In the light of Remark~\ref{r:51}, we prove as in the previous section:
\bprop
\lab{p:53}
For any $\varphi\in L^2(0,T; H^1_0(\Om))$, we have
\be
\lab{ccvm}
 \int_0^T\!\!\!\intb_\De u^\veps\varphi dxdt\to \int_\QT u\varphi dxdt.
\ee
\eprop

The homogenization result obtained in this case follows.

\bthm
\lab{c:liminf}
The limit  $u\in L^\infty(0,T;L^2(\Om))\cap L^2(0,T;H^1_0(\Om))$  of
(\ref{cvstari})--(\ref{cvgrii})  is the only solution of
\beq
\lab{pblim1inf}
&&
(1+a)\frac{\pa u}{\pa t}
-\Delta u  =
f\quad\mbox{in}\quad\QT,
\\
\lab{inituinf}
&&
u(0) = \frac{1}{(1+a)}u_0 + \frac{a}{(1+a)}v_0\quad\mbox{in}\quad\Om
\eeq
Moreover, the convergences in Proposition~\ref{p:uttauinf} hold on the whole
sequence and  $u\in C^0([0,T]; L^2(\Om))$, this being the sense of
(\ref{inituinf}).
\ethm
\bP
The proof of (\ref{pblim1inf}) is similar to the corresponding one  of the
Theorem~\ref{c:lim}. The test function
$\Phi^\veps$ is given by
$$
 \Phi^\veps = (1-w_\Rex)\varphi + w_\Rex G_\Ri
(\varphi),\quad\varphi\in\DD(\Om).
$$
The only interesting convergences are the following two:
$$
\left\vert\int_{\C_\veps^T}\na u^\veps (\na w_\Rex)(G_\Ri(\varphi) -
\varphi)dxdt\right\vert\leq C\vert \na u^\veps\vert_\QT \vert\na w_\Rex\vert_\QT
\left\vert G_\Ri(\varphi) - \varphi\right\vert_{L^\infty(\C_\veps^T)}\leq
$$
$$
 \leq C\gamma_\veps^{1/2}R_\veps =
C\left(\frac{r_\veps}{\veps}\right)^{1/2}\left(\frac{R_\veps}{\veps}\right)\to 0
$$

$$
 \left\vert\int_0^T\la f^\veps,w_\Rex (G_\Ri(\varphi) -
\varphi)\ra\right\vert\leq
C\left\vert (G_\Ri(\varphi) - \varphi)\na
w_\Rex\right\vert_{L^2(\C_\veps^T)}+C
\left\vert w_\Rex\na\varphi\right\vert_{L^2(\C_\veps^T\cup D_\veps^T)}\leq
$$
$$
 \leq C\left\vert\na\varphi\right\vert_{L^\infty(\Om)}\left(
\gamma_\veps^{1/2} R_\veps + \vert\C_\veps\cup\De\vert^{1/2}\right)\to 0,
$$
where we have used the a priori estimates of Proposition~\ref{p:23},
Remark~\ref{r:37} and  Proposition~\ref{p:web}.

Using Proposition~\ref{p:uttauinf} and hypotheses
(\ref{cvu0})--(\ref{cvv0}), we obtain with the same argument as before
$$
 \lim_{\veps\to 0}\int_\Om \rho^\veps u^\veps_0\Phi^\veps\eta(0) dx =
\eta(0)\int_\Om (u_0+ av_0)\varphi dx
$$
which achieves the proof.
\eP

\section{Homogenization in the case $r_\veps <<\veps^3$.}\lab{s:6}

As in this
case $\gamma_\veps\to 0$, we only can prove:

\bthm
\lab{p:uttau0f6}
There exists $ u\in L^\infty(0,T;L^2(\Om))\cap L^2(0,T;H^1_0(\Om))$  such that
\beq
\lab{cvstar0f6}
u^\veps  \buildrel{\star}\over\wto  u &\quad\mbox{in}\quad &
L^\infty(0,T;L^2(\Om))
\\
\lab{cvl20f6}
u^\veps  \wto  u & \quad\mbox{in}\quad & L^2(0,T;H^1_0(\Om))
\\
\lab{intbuef6}
\frac{1}{\vert\De\vert} u^\veps\chi_\De\to
v_0 &\quad\mbox{in} \quad &  \DD'(\Om)\quad\mbox{a.e.}\quad t\in
[0,T]
\eeq
where $u$ is the only solution of the following  problem:
\beq
\lab{pblim1inf06}
 \frac{\pa u}{\pa t}  - \Delta u & = & f\quad\mbox{in}\quad\QT
\\
\lab{inituinf06}
 u(0) & = & u_0\quad\mbox{in}\quad\Om
\eeq
\ethm
\bP
The convergences (\ref{cvstar0f6})--(\ref{cvl20f6}) hold on some subsequences;
they are insured by  Proposition~\ref{p:23}.  We have to remark that
(\ref{bornue}) is the hypothesis which insures the existence of
$v\in L^\infty(0,T;L^2(\Om))$ which satisfies
$$
 \frac{1}{\vert\De\vert} u^\veps\chi_\De\to
v \quad\mbox{in} \quad  \DD'(\Om)\quad\mbox{a.e.}\quad t\in
[0,T]
$$
(see Lemma~A-2
\cite{BellieudB}); we have to  prove that $v=v_0$.

Acting as usual, we take
\be
\lab{deffie}
\Phi^\veps = (1-w_{R_\veps})\varphi + w_{R_\veps} G_\Ri(\psi)
\ee
for some $R_\veps\in\RR_\veps$ and $\varphi,\psi\in\DD(\Om)$. Notice that in this
case we have
\be
\Phi^\veps\to \varphi\quad\mbox{in}\quad H^1_0(\Om)
\ee
because obviously $w_\Rex \to 0$ in $H^1_0(\Om)$.

Passing to the limit in the variational formulation, we obtain in a
straightforward manner
$$
 -\int_\QT u\varphi\eta' dxdt -a \int_\QT v\psi\eta' dxdt + \int_\QT\na
u\na\varphi\eta dxdt = \int_0^T\la f,\varphi\ra\eta dt +
$$
$$
+ \left(\int_\Om
u_0\varphi dx + a\int_\Om v_0\psi dx\right) \eta(0),\quad\forall\eta\in \DD([0,T[)
$$
Setting $\varphi = 0$, we find that $v$ is independent of $t$ and that $v\in
C^0([0,T]; L^2(\Om))$, which achieves $v = v_0$. Then, setting $\psi = 0$, we
prove (\ref{pblim1inf06}) and (\ref{inituinf06}), the last one holding also in the
sense of $C^0([0,T]; L^2(\Om))$.
\eP

\vs5
{\bf Acknowledgements.} This work was done during the visit of F.~Bentalha
and D.~Poli\c{s}evschi at the I.R.M.A.R.'s Department of Mechanics
(University of Rennes~1) whose support is gratefully acknowledged. Also,
this work corresponds to a part of the C.N.C.S.I.S. Research Program
33079-2004.



%
\newcommand{\noopsort}[1]{}

\vs5

*  University of Batna, Department of Mathematics, Batna, Algeria,

\vs5

** Universit\'e de Rennes1, I.R.M.A.R, Campus de Beaulieu,
35042 Rennes Cedex (France)

\vs5

*** I.M.A.R., P.O. Box 1-764, Bucharest (Romania).

\end{document}